\documentstyle[leqno,12pt]{amsart}

\newcommand{\half}{{\frac{1}{2}}}
\newtheorem{thm}{Theorem}

\title{A characterization of the Rogers q-Hermite polynomials}
       
\author{Waleed A. Al-Salam} 
\thanks{Research  was partially supported by NSERC (Canada) grant A2975} 
\address[Waleed. A. Al-Salam]{Department of Mathematics\\ University of Alberta\\
   Edmonton, Canada\\ T6G 2G1}
\email{waleed@@euler.math.ualberta.ca}

\begin{document}
\maketitle

\begin{abstract} 
In this paper we characterize the Rogers q-Hermite polynomials as the
only orthogonal polynomial set which is also ${\cal D}_q$-Appell
where ${\cal D}_q $ is the
Askey-Wilson finite difference operator.
\end{abstract}


\section{Introduction}
Appell polynomials sets $\{P_n(x)\}$ are generated by the  relation

\begin{equation}\label{gen-fun}
A(t) e^{xt} = \sum_{n=0}^\infty P_n(x)\; t^n,
\end{equation}
\noindent
where $A(t)$ is a formal power series in $t$ with $A(0)=1. $ This
definition implies the equivalent property that
\begin{equation}\label{diff-eq}
D P_n(x)= P_{n-1}(x), \hskip1in D=d/dx,
\end{equation}
Examples of such polynomial sets are

\begin{equation}\label{ex}
\left\{ \frac{x^n}{n!}\right\}, \left\{ \frac{B_n(x)}{n!}\right\}, 
  \left\{ \frac{H_n(x)}{2^n n!} \right\}
\end{equation}
where $B_n(x)$ is the nth Bernoulli polynomial and $H_n(x)$ is the $n$th
Hermite polynomials generated by

\begin{equation}
e^{2xt-t^2} = \sum_{n=0}^\infty H_n(x) \frac{ t^n}{n!} .
\end{equation}

By an orthogonal polynomial set (OPS) we shall mean those polynomial sets
which satisfy a three term recurrence relation of the form

\begin{equation}\label{3-term}
P_{n+1}(x)= (A_n x+B_n) P_n(x) - C_n P_{n-1}(x), \qquad (n=0,1,2,\cdots)
\end{equation}
with $P_0(x)=1,\ P_{-1}(x)=0,$ and $A_nA_{n-1}C_n > 0. $

By Favard's theorem \cite{chihara} this is equivalent to  the existence of a
positive measure $d\alpha(x)$ such that
\begin{equation}
\int_{-\infty}^\infty P_n(x) P_m(x)\; d\alpha(x) = K_n \delta_{nm}.
\end{equation}
As we see from the examples \eqref{ex} some Appell polynomials are
orthogonal and some are not. This prompted Angelesco \cite{angelesco}
to prove that
{\it the
only orthogonal polynomial sets which are also Appell is the Hermite
polynomial set.}
This theorem was rediscovered by several authors later on
(see, e.g., \cite{shohat}).

There were several extensions and/or analogs of Appell polynomials that
were introduced later. Some are based on changing the operator $D$ in
\eqref{diff-eq} into another differentiation-like operator or by replacing the
generating relation \eqref{gen-fun} by a more general one. In most of these
cases theorems like Angelesco's were given. For example Carlitz
\cite{carlitz}
proved that the Charlier polynomials are the  only OPS which satisfy
 the difference relation
\begin{equation}
 \Delta P_n(x) = P_{n-1}(x), \qquad\qquad (\Delta f(x)=f(x+1)-f(x).)
\end{equation}
See \cite{wal} for many other references.

A new and very interesting analog of Appell polynomials were introduced
recently, as a biproduct of other considerations, by Ismail and Zhang 
\cite{ism-zhang}. 
In discussing the Askey-Wilson operator they defined a new q-analog of the
exponential function $e^{xt}.$ This we describe in the next section.


\section{Notations and Definitions}
\setcounter{equation}{0}

The Askey-Wilson operator is defined by
\begin{equation}
{\cal D}_q f(x)= \frac{\delta_q f(x)}{\delta_q x},
\end{equation}
where $x=\cos\theta$ and
\begin{equation}
 \delta_q g(e^{i\theta})= g(q^{1/2} e^{i\theta})- g(q^{-1/2}e^{i\theta}).
\end{equation}
We further assume that $-1<q<1$ and use the notation

\begin{eqnarray}
(a;q)_0 &=& 1,\qquad (a;q)_n=(1-a)(1-qa)\cdots (1-aq^{n-1}),
   \quad (n=1,2,..) \\
(a;q)_\infty &=& \prod_{k=0}^\infty (1-aq^j).
\end{eqnarray}

There are two q-analogs of the exponential function $e^x$ given by the infinite
products
\begin{equation}\label{qexp-1}
e_q(x) = \frac{1}{(x;q)_\infty}=\sum_{k=0}^\infty \frac{x^k}{(q;q)_k},
\end{equation}
and
\begin{equation}\label{qexp-2}
\frac{1}{e_q(x)} = (x;q)_\infty = \sum_{k=0}^\infty (-1)^k
   \frac{q^{\half k(k-1)}}{(q;q)_k } x^k.
\end{equation}

We shall also use the function
\begin{equation}
\Psi_n(x) = i^n (iq^{(1-n)/2}e^{i\theta};q)_n (iq^{(1-n)/2}
e^{-i\theta};q)_n,
\end{equation}
so that 
\begin{eqnarray}
\Psi_{2n}(x) &=&  \prod_{k=0}^{n-1} \left[ 4x^2+(1-q^{2n-1-2k})(1-q^{
1-2n+2k})  \right] \nonumber \\
\Psi_{2n+1}(x) &=& 2x  \prod_{k=0}^{n-1} \left[ 4 x^2 - (1-q^{2n-2k})
   (1-q^{-2n+2k}) \right] \nonumber \\
4 x^2 \Psi_n(x) &=& \Psi_{n+2}(x)+ (1-q^{n+1})(1-q^{-n-1})\Psi_n(x)
\end{eqnarray}

Thus
\begin{equation}\label{psi}
{\cal D}_q \Psi_n(x)= 2 q^{(1-n)/2} \frac{1-q^n}{1-q}\ \Psi_{n-1}(x).
\end{equation}
and
\begin{equation}\label{psi-1}
{\cal D}_q\left[ x\; \Psi_n(x)\right] = \frac{q^{(1+n)/2} -q^{-(n+1)/2}}{
   q^{1/2}-q^{-1/2}}\ 2x\; \Psi_{n-1}(x).
\end{equation}

Iterating \eqref{psi} we get
\begin{equation}\label{psi-k}
{\cal D}_q^k \Psi_n(x)=\ 2^k q^{\frac14 k(k+1)-\half nk} \frac{(q;q)_n}{
   (q;q)_{n-k} (1-q)^k} \Psi_{n-k}(x).
\end{equation}

The Ismail-Zhang q-analog of the exponential function \cite{ism-zhang} is
\begin{equation}\label{Dq-Appell}
{\cal E}(x) =\sum_{n=0}^\infty \frac{q^{n(n-1)/4} (1-q)^n}{2^n (q;q)_n}
   \ \Psi_n(x)\; t^n.
\end{equation}

It follows from \eqref{Dq-Appell} and \eqref{psi} that
\begin{equation}
 {\cal D}_q {\cal E}(x) = t\; {\cal E}(x).
\end{equation}

This suggested to Ismail and Zhang to define the $\cal D_q$-Appell
polynomials as those, in analogy with \eqref{gen-fun}, defined by
\begin{equation}\label{gen-fun-2}
 A(t){\cal E}(x) = \sum_{n=0}^\infty P_n(x)\ t^n ,
\end{equation}
so that
\begin{equation}\label{diff-rel}
   {\cal D}_q P_n(x)= P_{n-1}(x).
\end{equation}
An example of such a set is the Rogers q-Hermite polynomials,
$\{H_n(x|q)\}$,  (see \cite{wal-ism,ask-ism,ism-stan}).
\begin{equation}\label{gen-fun-H}
\prod_{n=0}^\infty \left(1-2xt q^n+t^2 q^{2n}\right)^{-1}=
      \sum_{n=0}^\infty H_n(x|q) \frac{t^n}{(q;q)_n}.
\end{equation}

They satisfy the three term recurrence relation
\begin{equation}\label{rec-rel-H}
H_{n+1}(x|q)= 2x H_n(x|q)-(1-q^n)H_{n-1}(x|q),\qquad n=0,1,2,3,...
\end{equation}
with $H_0(x|q)=1,\ H_{-1}(x|q)=0.$

\vskip0.4in

\section{The Main Result}
\setcounter{equation}{0}
We now state our main result:
\begin{thm}
The orthogonal polynomial sets which are also ${\cal D}_q$-Appell, i.e.,
satisfy \eqref{diff-rel} or \eqref{gen-fun-2} is the set of the Rogers q-Hermite polynomials.
\end{thm}

{\bf Proof} Let $\{Q_n(x)\}$ be a polynomial set which is both orthogonal
and ${\cal D}_q$-Appell. That is $\{Q_n(x)\}$ satisfy \eqref{gen-fun-2}
and \eqref{3-term}.

We next note that \eqref{gen-fun-H} implies that
\begin{equation}
h_n(x|q) = \frac{(1-q)^n q^{n(n-1)/4}}{2^n (q;q)_n} H_n(x|q)
\end{equation}
satisfy

\begin{equation}
  {\cal D}_q h_n(x|q)= h_{n-1}(x|q),
\end{equation}
so that $\{h_n(x|q)\}$ is a ${\cal D}_q-$Appell polynomial set and at the
same time is an OPS satisfying the three term recurrence relation

\begin{equation}\label{3-term-h}
(1-q^{n+1}) h_{n+1}(x|q)= (1-q)q^{n/2} x h_n(x|q)-\frac14 (1-q)^2
q^{n-1/2} h_{n-1}(x|q)
\end{equation}

It also follows from \eqref{gen-fun-2} that any two polynomial sets 
$\{R_n(x)\}$ and $\{ S_n(x)\}$, in
that class are related by $R_n(x)=\sum_{k=0}^n c_{n-k} S_k(x).$ Thus the
solution to our problem may be expressed as

\begin{equation}\label{Q}
 Q_n(x)= \sum_{k=0}^n a_{n-k} h_k(x|q).
\end{equation}
for some sequence of real constants $\{a_n\}.$ We may assume without loss
of generality that $a_0=1.$

The three term recurrence relation satisfied by $\{Q_n(x)\}$ is

\begin{equation}\label{3-term-Q}
(1-q^{n+1})Q_{n+1}(x)= \left( (1-q)q^{n/2} x + \beta_n\right)Q_n(x)
  -\gamma_n Q_{n-1}(x),
\end{equation}
with $Q_0(x)=1,\ Q_{-1}(x)=0.$ Thus $Q_1(x)=x+\beta_0=a_1+h_1(x|q),$ from
which it follows that $a_1=\beta_0.$

Putting \eqref{Q} in \eqref{3-term-Q} and using \eqref{3-term-h} to
replace $xh_k(x|q)$ in terms of $h_{k+1}(x|q)$ and $h_{k-1}(x|q)$ we get,
on equating coefficients of $h_k(x|q),$

\begin{equation}\label{sys-eq}
(1-q^{(n-k+1)/2})(1+q^{(n+1+k)/2}) a_{n+1-k} -\beta_n a_{n-k}
   +\left[ \gamma_n -\frac14(1-q)^2 q^{(n+k)/2} \right] a_{n-k-1}=0,
\end{equation}
valid for all $n$ and $k=0,1,2,...,n+1$ provided we interpret
$a_{-1}=a_{-2}=0.$ It is easy to see that this system of equations is
equivalent to the solution of our problem.

Putting $k=n$ in \eqref{sys-eq} we get
\begin{equation}\label{beta}
        \beta_n = (1-q^\half)(1+q^{n+\half})a_1.
\end{equation}
Hence if $\beta_0=0$ then $\beta_n=0$ for all $n.$ In fact if $\beta_m=0$
for any $n=m$ then $\beta_n=0 $ for all $n.$

Now we treat these two cases seperately.

{\bf Case I.} ($\beta_0=0$).

The system \eqref{sys-eq} can now be written as

\begin{equation}\label{sys-2}
(1-q^{(k+1)/2}) (1+q^{n+\half (1-k)}) a_{k+1}+
   \left[ \gamma_n -\frac14 (1-q)^2 q^{n-\half k}\right] a_{k-1}=0.
\end{equation}

Since $a_1=0 $ then it follows from \eqref{sys-2} that
$a_{2k+1}=0$ for all $k.$ In particular we get

\begin{equation}
   \gamma_n = \frac14 (1-q)^2 q^{n-\half} - a_2 (1-q)(1+q^n),
\end{equation}
so that if $a_2=0$ then 
\begin{equation}
Q_n(x)= h_n(x|q).
\end{equation}

Now we show that $a_2\ne 0$ leads to contradiction. To do this replace
$k$ by $2k-1$.  We get

\begin{equation}\label{eq-3}
(1-q^k)(1+q^{n-k+1})a_{2k}+\left[ \frac14 (1-q)^2 q^{n-\half} (1-q^{1-k})
  -a_2 (1-q)(1+q^n) \right] a_{2k-2} = 0.
\end{equation}

Keep $k$ fixed and let $n\to\infty.$  We get
$(1-q^k)a_{2k}=(1-q)a_2 a_{2k-2}.$ Thus

\begin{equation}
  a_{2k} = \frac{(1-q)^k}{(q;q)_k} a_2^k.
\end{equation}

Putting this value in \eqref{eq-3} we get $q^{1-k}=1.$ This is a
contradiction and Case I is finished.

{\bf Case II} ($\beta_0\ne 0$).

We start with \eqref{sys-eq} we get, assuming $a_1\ne 0$, 
\begin{equation}
\gamma_n= \frac14 (1-q)^2 q^{n-\half}+(1-q^\half)(1+q^{n+\half})a_1^2
   - (1-q)(1+q^n)a_2.
\end{equation}
Putting this value of $\gamma_n$ and the value of $\beta_n$ in
\eqref{beta} in \eqref{sys-eq}, and finally equating coefficients of $q^n$
and the terms independent of $n$ we get the pair of equation systems
\begin{equation}\label{eq-6}
(1-q^{(k+1)/2})a_{k+1}-(1-q^\half)a_1a_k+\left\{
(1-q^\half)a_1^2-(1-q)a_2\right\} a_{k-1} =0 
\end{equation}
and
\begin{eqnarray}\label{eq-5}
\lefteqn{ (1-q^{(k+1)/2})a_{k+1}-(1-q^\half)q^{k/2} a_1a_k +} \\
 & & \left\{ \frac14 (1-q)^2
q^{-\half}(q^{(k-1)/2}-1)+q^{k/2}(1-q^\half ) a_1^2
      -(1-q)q^{(k-1)/2}a_2  \right\} a_{k-1} =0 \nonumber
\end{eqnarray}

Eliminating $a_{k+1}$ in these equations we get

\begin{eqnarray}\label{eq-4}
\lefteqn{
   (1-q^\half)(1-q^{k/2}) a_1a_k + \left\{ (1-q) a_2 (1-q^{(k-1)/2})
       \right. }\\
& &  - \left. (1-q^\half)(1-q^{k/2}) a_1^2 -\frac14 (1-q)^2 q^{-\half}
(1-q^{(k-1)/2} )
  \right\} a_{k-1} = 0. \nonumber
\end{eqnarray}
This equation is of the form $(1-q^{k/2})a_1a_k= c(1-b q^{k/2})a_{k-1}$ so
that the general solution of \eqref{eq-4} is

\begin{equation}
a_k = c^k \frac{(bq^\half;q^\half)_k}{(q^\half;q^\half)_k}
\end{equation}

Putting this in \eqref{eq-6} we get that $b=0.$ On the other hand 
\eqref{eq-5} gives
that $c^2= \frac14 (1-q)^2 q^{-\half}. $ Finally putting those values of
$a_k$ in (3.13) we get that $\gamma_n =0$ which is a contradiction.

This completes the proof of the theorem.


\section{Generating Function}
\setcounter{equation}{0}

We obtain, for the q-Hermite polynomials, a generating function of the
form \eqref{gen-fun-2}. More specifically we prove
\begin{thm}
Let $H_n(x|q)$ be the $n$th Rogers q-Hermite polynomial. Then we have
\begin{equation}\label{thm-2}
\sum_{n=0}^\infty \frac{q^{n(n-1)/4}}{(q;q)_n}\; H_n(x|q)\; t^n
  = (t^2 q^{-\half};q^2)_\infty {\cal E}(x).
\end{equation}
\end{thm}

{\bf Proof.} Let $A(t)=1+a_1 t+a_2 t^2+a_3 t^3+\cdots $ and
\begin{equation}\label{gen-fun-h}
  A(t){\cal E}(x) = \sum_{n=0}^\infty h_n(x|q) t^n.
\end{equation}

Then we get
\begin{equation}\label{h}
  h_n(x|q)=\sum_{k=0}^n a_{n-k} c_k \Psi_k(x).
\end{equation}
where

\begin{equation}
      c_k = \frac{(1-q)^k}{2^k (q;q)_k}\; q^{k(k-1)/4}.
\end{equation}

To calculate the coefficients $\{a_n\}$ we first iterate \eqref{3-term-h}
we get
\begin{eqnarray}\label{rec-h}
4x^2 h_n(x|q) &=& \frac{4}{(1-q)^2} (1-q^{n+1})(1-q^{n+2})q^{-n-\half}
  h_{n+2}(x|q) \\
  &+& (2-q^n-q^{n+1})h_n(x|q)+ \frac{(1-q)^2}{4} q^{n-\frac32}
     h_{n-2}(x|q). \nonumber
\end{eqnarray}
 Putting \eqref{h} in \eqref{rec-h}, using (2.6) and then equating coefficients of
$\Psi_k(x)$ we get after some simplification
\begin{eqnarray}\label{rec-3}
 \lefteqn{\frac{4}{(1-q)^2} q^{-n-\half} (1-q^{n-k+2})(1-q^{n+k+1})
a_{n+2-k}+} \\
  & & q^{-k-1} \left\{ 1+q^{2k+2}-q^{n+k+1}-q^{n+k+2}\right\} a_{n-k}+
       \nonumber \\
 & & \frac{(1-q)^2}{4} q^{n-\frac32} a_{n-2-k}=0 \qquad (k=0,1,...,n+2).
      \nonumber
\end{eqnarray}

By direct calculation of $a_1,\; a_2,\; a_3$ we see easily that
$a_1=a_3=0.$ Thus \eqref{rec-3} shows that $a_{2k+1}=0$ for all $k.$

Furthermore we can easily verify that
\begin{equation}
a_{2j} = (-1)^j \frac{(1-q)^{2j}}{2^{2j}(q^2;q^2)_j}\  q^{j(j-\frac32)}
 \qquad (j=0,1,2,3,...
\end{equation}

Hence

\begin{eqnarray}
A(t) &=& \sum_{j=0}^\infty (-1)^j \frac{q^{j(j-1)}}{(q^2;q^2)_j}
       \left( \frac{(1-q)^2 t^2}{4} q^{-\half}\right)^j\\
  &=& \left( \frac{(1-q)^2}{4} t^2 q^{-\half};q^2\right)_\infty.\nonumber
\end{eqnarray}
After some rescaling we get the theorem.

As a corollary of \eqref{thm-2} we state the pair of inverse relations 
\begin{eqnarray}
\Psi_n(x) &=& \sum_k \frac{(q;q)_n q^{k(k-n)}}{(q^2;q^2)_k (q;q)_{n-2k}}
     H_{n-2k}(x|q), \\
H_n(x|q) &=& \sum_k (-1)^k \frac{(q;q)_n\ q^{k(2k-n-1)}}{
   (q^2;q^2)_k(q;q)_{n-2k}} \Psi_{n-2k}(x). \label{hn}
\end{eqnarray}
These follows from the identities \eqref{qexp-1} and \eqref{qexp-2}

Formula \eqref{hn} and \eqref{psi-k} give
\begin{equation}
H_n(x|q)\ =\ \frac{1}{e_{q^2}(\frac{(1-q)^2}{4}q^{-\half}{\cal D}_q^2)}
  \Psi_n(x).
\end{equation}
This is a q-analog of the formula
$$ e^{-D^2} x^n = H_n(x) $$
for the regular Hermite polynomials (1.4).


\end{document}